\documentclass[12pt]{article}
\usepackage{amssymb}
\usepackage[usenames]{color}
\usepackage[colorlinks=true,
linkcolor=webgreen,
filecolor=webbrown,
citecolor=webgreen]{hyperref}

\definecolor{webgreen}{rgb}{0,.5,0}
\definecolor{webbrown}{rgb}{.6,0,0}

\usepackage{epsfig}
\usepackage{multirow}

\begin{document} 

\begin{center}
\vspace*{0.1in}
{\huge Solitaire: Recent Developments}
\vskip 20pt
{\bf John D. Beasley}\\
September 2003\footnote{Original version at
\href{http://gpj.connectfree.co.uk/gpjj.htm}
{\tt http://gpj.connectfree.co.uk/gpjj.htm} \\
\hspace*{0.25in}Converted to \LaTeX ~by George Bell with minor
modifications to the text, November 2008.} \\
{\tt johnbeasley@mail.com}\\
\end{center}
\vskip 30pt 
\centerline{\bf Abstract}
\noindent
This special issue on Peg Solitaire has been put
together by John Beasley as guest editor,
and reports work by John Harris, Alain Maye, Jean-Charles Meyrignac,
George Bell, and others.
Topics include: short solutions on the $6\times 6$ board
and the 37-hole ``French" board,
solving generalized cross boards and long-arm boards.
Five new problems are given for readers to solve,
with solutions provided.

\pagestyle{myheadings}
\markright{The Games and Puzzles Journal---Issue 28, September 2003\hfill}

\thispagestyle{empty} 
\baselineskip=15pt 
\vskip 30pt 

\section{Introduction and historical update}

There has recently been a flurry of activity on the game of Peg Solitaire,
and I have suggested to George Jelliss that
\href{http://www.gpj.connectfree.co.uk/index.htm}{The Games and Puzzles Journal} \cite{GPJ}
might be a convenient place for people to report new discoveries.
His reaction was that he would like to introduce the game to readers by
dedicating a special number to it, after which he will consider
contributions as they arise, and he has asked me to provide the material
for this special edition.
It updates the material given in my book
\textit{The Ins and Outs of Peg Solitaire} \cite{Beasley}
and what was given there will not normally be repeated here,
but enough background will be given to put any reader not previously
familiar with the game's development in the picture.
\textit{The Ins and Outs} is now out of print and will probably remain so,
but it can be found in most academic and many UK public service libraries,
and there appears to be a steady trickle of copies on the secondhand
market\footnote{Try \href{http://www.ABEbooks.com}{ABEbooks.com}.}.
The 1992 edition differs from the 1985 only in the addition of a page
summarizing intervening developments and discoveries,
and I can supply photocopies of this on request.

\noindent
The game's historical background is now well known.
It originated in France in the late seventeenth century
(there are references in French sources back to 1697),
and it appears to have been the ``Rubik's Cube" of the court of Louis XIV.
I summarized its early history \cite[p. 3--7]{Beasley}
and little appears to have
been discovered since, but one statement in the book now needs modification.
I took a very cautious view of a passing reference in a letter written by
Horace Walpole in 1746, fearing that it might have referred to a card game,
but David Parlett, who has looked into the games of the period much more
deeply and extensively than I, tells me that my fears were groundless:
``Patience dates from the late eighteenth century, did not reach England
until the nineteenth, and was not called Solitaire when it did"
\cite[p. 157]{Parlett}.
So the spread of our Solitaire to England by the middle of the
eighteenth century can be taken as established.

\noindent
There is one matter in which discovery remains conspicuous by its absence.
It has frequently been written that the game was invented by a prisoner
in the Bastille, but I reported in 1985 that the earliest reference
to this appeared to be in an English book of 1801, and nobody has yet
drawn my attention to anything earlier.
An uncorroborated English source
of 1801 is of course quite valueless as evidence for an alleged occurrence
in France over a century before, and anyone who repeats this tale without
citing a French source earlier than 1801 should regard himself as
perpetuating myth rather than history.
Sadly, the more picturesque a legend surrounding the origin of a game
or puzzle, the greater the likelihood that somebody has invented
it along the way.

\section{The 6 x 6 board: the work of John Harris}

I gave solutions to various problems on the $6\times 6$ board in the
1985 edition of \textit{The Ins and Outs}, and on page 252 of the 1992
edition I added a note that John Harris had found all possible
15-move solutions, one by hand and the rest by computer.
I refrained from giving details on the grounds that he might still
wish to publish them himself, but to the best of my knowledge
he has not done so, and others are beginning to reproduce his results.
I therefore think I should summarize what he sent me in 1985--86,
if only to establish his priority.

\begin{figure}[htb]
\centering
\epsfig{file=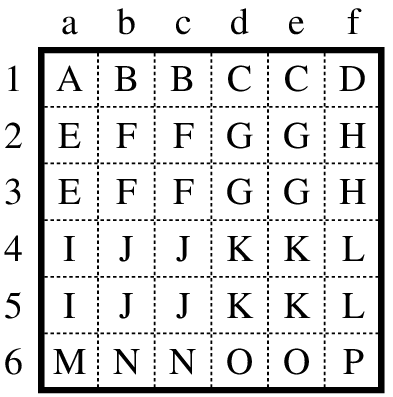}
\caption{The $6\times 6$ square board divided into 16 ``Merson regions" A--P.}
\label{fig1}
\end{figure}

\noindent
Robin Merson observed back in 1962 that the holes of the $6\times 6$
board could be divided into 16 regions such that only the first jump
of a multi-jump move could open up a new region (Figure~\ref{fig1});
any later jump had to be between regions already opened.
It follows that it takes at least 15 moves to clear the board if
the initial vacancy is in a non-corner square,
and 16 moves if it is in a corner (because the first move refills
this corner and we are still left with 15 regions to be opened).
Harris found a 16-move solution to the problem
``vacate a1 and play to finish there" back in 1962,
and Harry O. Davis subsequently found 15-move solutions
to the problems ``vacate c1, finish at f4"; and
``vacate c3, finish at f6".
All these are in \textit{The Ins and Outs}.
That to ``start and finish at a1" ends with an elegant eight-sweep loop.

\noindent
Subsequently (letter to me dated 26 August 1985) Harris found a 15-move
solution to the problem ``vacate c3, finish at c6":
c1-c3, a1-c1, d1-b1, f1-d1, a3-a1-c1-e1 (5), a5-a3, c4-a4-a2-c2-c4,
d4-b4, c6-c4-a4, e3-c3-a3-a5-c5 (10), f3-f1-d1-d3, f5-f3,
d6-d4-d2-f2-f4-d4, f6-d6, a6-c6-e6-e4-c4-c6.
``Don't know how to find these," he wrote,
``just copied a Davis beginning and got lucky."
Harris then attacked the problem by computer, and by August 1986
he had found 15-move solutions to all the problems with non-corner
starts apart from ``start and finish at c1".
His computer proved this to require 16 moves.
The remaining corner-start problems, ``vacate a1, finish at a4 or d4"
had been solved in 16 moves by Davis, and his solutions had appeared
in the instructions to Wade Philpott's 1974 game SWEEP.
Harris's solution to ``start and finish at c3" used 5 single jumps,
then 6 double jumps, then 4 jumps from corners:
c1-c3, a2-c2, d2-b2, d4-d2, d6-d4 (5), f3-d3-d5, b6-d6-d4,
f5-d5-d3, e1-e3-e5, c4-c2-e2, a4-c4-c6 (11), a6-a4-a2-c2,
f6-d6-b6-b4-b2-d2-d4, f1-f3-f5-d5-d3, a1-c1-e1-e3-c3.
``This has to be my favorite solution," he wrote.

\noindent
Readers who revel in the power of modern computers may care to note
that all this was done on a TRS-80, which if memory serves me right
offered a mere 64Kb of RAM for operating system, program,
and data together, backed up by a single 52Kb disc drive.
Harris's results have recently been confirmed by
Jean-Charles Meyrignac.

\section{Solutions on the classical 33-hole and 37-hole boards}

The classical 33-hole (3-3-7-7-7-3-3) and 37-hole (3-5-7-7-7-5-3) boards
offer no simple test for minimality such as is provided by the need to
open up each of Merson's regions on a $6\times 6$ board,
and there is usually a gap of two or three between the length
of the shortest solution actually discovered and the number
of moves that can be proved necessary by simple means.
This gap can be filled only by an exhaustive analysis by computer.

\subsection{The 33-hole board}

In the original 1985 edition of \textit{The Ins and Outs},
I listed the shortest solutions found by Ernest Bergholt
and Harry O.
Davis to the single-vacancy single-survivor
problems on the standard 33-hole board,
and I reported some last-minute computer calculations by myself
which demonstrated them indeed to be the shortest possible.
However, this was ``proof of non-existence by failure to find 
despite a search believed exhaustive",
and to achieve it on the machine at my disposal I had to resort
to some fairly complicated testing to identify and reject blind alleys.
I therefore took the view that the proof should be regarded as
provisional pending independent confirmation.

\noindent
No such confirmation had been reported to me when the 1992
edition went to press, but on 24 October 2002 Jean-Charles Meyrignac
reported \cite{Meyrignac}
that he had programmed the calculation independently and had
verified that the solutions of Bergholt and Davis were indeed
optimal\footnote{The report on his web site merely said ``All solutions",
but he has clarified the matter in an e-mail to me.}.
My 1984 machine offered only 32Kb of RAM for program and data together,
even less than that provided by Harris's TRS-80,
though I did have two 100Kb disc drives.
Meyrignac, with a more powerful present-day machine at his disposal,
had no need for complicated restriction testing and could perform
a complete enumeration, reproducing all known solutions as well
as demonstrating that there were none shorter.

\subsection{The 37-hole ``French" board (see Figure~\ref{fig5}a)}

Although this was historically the first board to be used,
minimal solutions on it appear to have received less attention
than those on the 33-hole board, and in \textit{The Ins and Outs}
I could only report some relatively recent findings by
Leonard Gordon and Harry O. Davis.
Four of Gordon's solutions were subsequently beaten by Alain Maye
(work dating from 1985-86 but only recently brought to my notice),
and I would have reported this in the 1992 edition of
\textit{The Ins and Outs} had I been aware of it.
Meyrignac has now performed an exhaustive enumeration by computer,
which shows that the problem ``vacate c1, play to a single survivor"
can be solved in 20 moves irrespective of which of the holes
b4/e1/e4/e7 is chosen to receive the survivor (Gordon and Maye
had got each case down to 21), and proves the remaining solutions
of Gordon, Davis, and Maye to be optimal.
Table~\ref{table1} below has been proved by Meyrignac to be definitive.

\begin{table}[htbp]
\begin{center} 
\begin{tabular}{ | c | c | c | l | }
\hline
Vacate & Finish & Length & Investigator \\		
\hline
\hline
\multirow{4}{*}{c1} & e1 & 20 & \multirow{4}{*}{Meyrignac (by computer)} \\
 & b4 & 20 & \\
 & e4 & 20 & \\
 & e7 & 20 & \\
\hline
\multirow{3}{*}{d3} & d2 & 21 & Gordon \\
 & a5 & 21 & Davis \\
 & d5 & 21 & Maye \\
\hline
\multirow{3}{*}{d6} & d2 & 20 & Gordon \\
 & a5 & 20 & Gordon \\
 & d5 & 20 & Maye \\
\hline
\end{tabular}
\caption{Summary of shortest solutions on the 37-hole French board.} 
\label{table1}
\end{center} 
\end{table}

\noindent
Maye's solutions:\newline
Vacate d3, finish at d5: d1-d3, b2-d2, d3-d1, f2-d2, e4-e2 (5), c4-c2,
a3-c3, d1-d3-b3, g3-e3, a5-a3-c3 (10), b5-b3-d3-f3, g5-g3-e3, d5-b5,
b6-b4, c7-c5 (15), c1-c3, f5-f3-d3-b3-b5-d5-f5, e1-e3, f6-f4, e7-e5 (20),
d7-d5-d3-f3-f5-d5.\newline
Vacate d6, finish at d5: b6-d6, c4-c6, c7-c5, a4-c4-c6, e7-c7-c5 (5),
e6-c6-c4, b2-b4, d3-b3, c1-c3-c5, a3-c3 (10), e4-e6, f6-d6, g5-e5,
e2-c2-c4-c6-e6-e4-e2, d2-f2 (15), g3-e3, g4-e4-e2-c2, a5-c5, e1-c1-c3,
d5-b5-b3-d3-d5.\newline
The most interesting of Meyrignac's is ``vacate c1, finish at e4",
which ends with an eight-sweep: e1-c1, d3-d1, b3-d3, c5-c3, c7-c5 (5),
e4-c4-c6, f2-d2-d4, b2-d2, g4-e4-e2, g3-e3 (10), a5-c5, f6-f4,
d5-f5-f3-d3-b3, c1-e1-e3, a3-c3 (15), e7-e5, d7-d5-f5, g5-e5,
b6-d6, a4-c4-e4-e2-c2-c4-c6-e6-e4!\newline
All Maye's and Meyrignac's solutions can be found on Meyrignac's
web site \cite{Meyrignac}.

\section{Generalized cross boards and long-arm boards}

\subsection{Generalized cross boards}

\begin{figure}[htb]
\centering
\epsfig{file=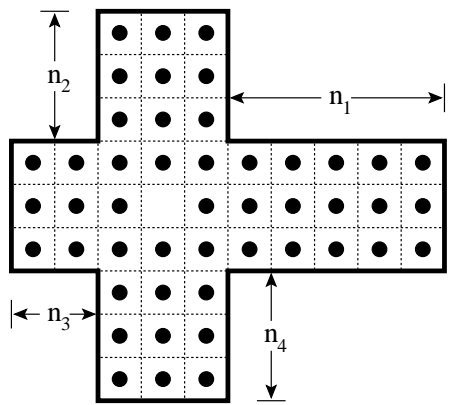}
\caption{The 48-hole generalized cross board $n_1=5$, $n_2=3$, $n_3=2$, $n_4=3$.}
\label{fig2}
\end{figure}

George Bell has been studying a class of boards he calls
``generalized cross boards".
These have a similar cross shape to the standard 33-hole board,
but the four $3\times n$ ``arms" are allowed to have different lengths
$n_1, n_2, n_3, n_4$ (including zero).
The standard 33-hole board is of course such a board ($n_1=n_2=n_3=n_4=2$),
as is Wiegleb's 45-hole board ($n_1=n_2=n_3=n_4=3$).
Shown in Figure~\ref{fig2} is a 48-hole example with
$n_1=5$, $n_2=3$, $n_3=2$, $n_4=3$.

\noindent
All these generalized cross boards are built up from rows of three,
so they are automatically null-class boards.
We can therefore hope that ``single-vacancy complement problems",
where we play to leave a single peg in the hole initially vacated,
will be solvable, and we shall describe a board as``solvable at X" if the
problem ``vacate X, play to finish at X" is solvable on it.
Making extensive use of the computer for investigation,
George has shown that there are exactly 12 generalized cross boards
which are solvable at every location.
Table~\ref{table2} lists all such boards---they range in size
from 24 to 42 holes,
and of course they include the standard 33-hole board
(but not Wiegleb's board, which is not solvable at the
middle square at the end of an arm).

\begin{table}[htbp]
\begin{center} 
\begin{tabular}{ | c | c | c | c | c | l | l |}
\hline
$n_1$ & $n_2$ & $n_3$ & $n_4$ & Holes & Symmetry & Comment \\		
\hline
\hline
2 & 1 & 2 & 0 & 24 & Lateral & \\
\hline
2 & 1 & 2 & 1 & 27 & Rectangular & \\
2 & 2 & 1 & 1 & 27 & Diagonal & \\
3 & 2 & 0 & 1 & 27 & & \\
\hline
3 & 2 & 1 & 1 & 30 & & \\
\hline
2 & 2 & 2 & 2 & 33 & Square & The standard 33-hole board \\
3 & 2 & 2 & 1 & 33 & & \\
\hline
3 & 3 & 2 & 1 & 36 & & \\
3 & 2 & 3 & 1 & 36 & Lateral & \\
\hline
3 & 2 & 3 & 2 & 39 & Rectangular & ``semi-Wiegleb" \\
3 & 3 & 2 & 2 & 39 & Diagonal & \\
\hline
3 & 3 & 3 & 2 & 42 & Lateral & \\
\hline
\end{tabular}
\caption{The 12 generalized cross board solvable at every location.} 
\label{table2}
\end{center} 
\end{table}

\noindent
Most of these problems are easy, but some are not.
Perhaps the hardest is given by the middle square at the end of a
long arm on the 39-hole board ``3,2,3,2",
which has two ``standard" arms and two ``Wiegleb" arms.
This is presented as a \textit{problem to solve} in the last section,
and its solution is unique to within symmetry and order
of jumps \cite{BellBeasley}.

\noindent
No generalized cross board other than these twelve is solvable everywhere.
George demonstrates this by applying conventional analysis to show that
no such board with an arm of length 5 or more can be solvable everywhere
(via the same technique as he uses for the general 6-arm case below),
and then performing a relatively simple and quick computer analysis of
the 45 remaining cases.
However, the computer analysis must be laboriously run over each case individually,
and he stresses that the results await independent verification.
His analysis of Wiegleb's board confirms my own \cite[p. 199--201]{Beasley}.

\subsection{Boards with longer arms}

\begin{figure}[htb]
\centering
\epsfig{file=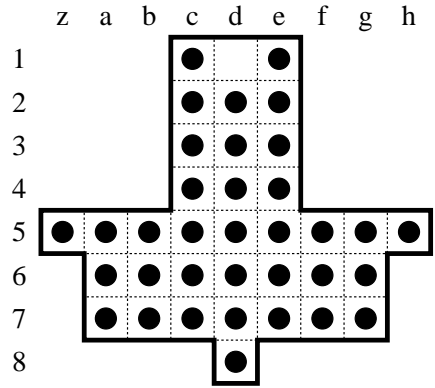}
\caption{The 36-hole ``mushroom board" with a solvable d1-complement.}
\label{fig3}
\end{figure}

The investigation above showed that no generalized cross board with an
arm longer than three was solvable everywhere,
but George wondered what would happen if a longer arm were attached to
a board of some other shape.
He came up with the 36-hole ``mushroom board" (Figure~\ref{fig3}),
which proved that a board with a 4-arm could be solvable at all locations,
in particular at the middle of the end of the arm
(always likely to be the most difficult square).
For convenience, we invert the mushroom so that this key square is at the top,
and we continue to call it ``d1", adding a z-file to the left of the a-file.

\noindent
The d1-complement on this board can be solved by
d3-d1, d5-d3, b5-d5, d6-d4-d2, f5-d5, d8-d6-d4, e3-e5, e1-e3, e6-e4-e2,
h5-f5, g7-g5-e5, b7-d7, e7-c7-c5, c4-c6, b6-d6, z5-b5, a7-a5-c5,
c2-c4-c6-e6-e4, c1-e1-e3-e5, f7-f5-d5-d3-d1.
The a6-complement is another tricky one
(it fails if the arm is only of length 2, or is absent altogether),
but the other single-vacancy complement problems are not difficult.

\noindent
This board has only lateral symmetry, but it is not difficult to
construct 4-arm boards solvable everywhere that have square symmetry.
One example is the 129-hole board obtained by taking a $9\times 9$ square
and attaching a 4-arm to the middle of each side.

\begin{figure}[htb]
\centering
\epsfig{file=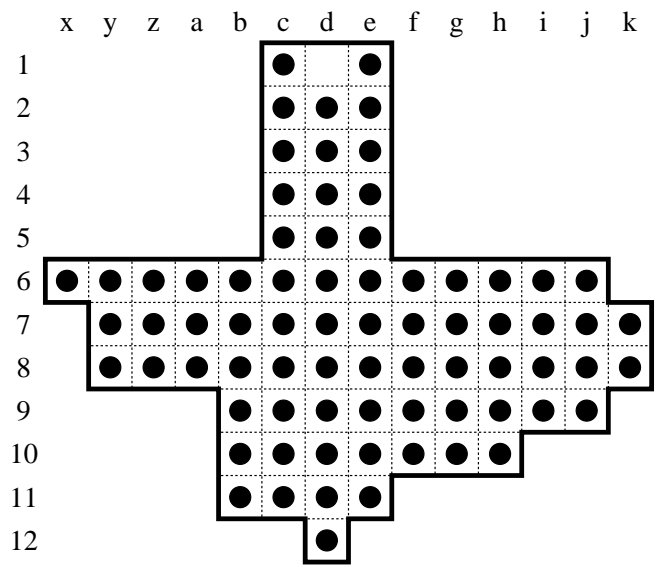}
\caption{A 75-hole board with a solvable d1-complement.}
\label{fig4}
\end{figure}

\noindent
Initial experimentation suggested that any board with a 5-arm would
be unsolvable at the mid-end of the arm, but a proof covering all cases
was elusive and eventually we found a 90-hole board which was solvable there.
Subsequent exploration brought the number of holes down to 75,
and further reduction may be possible.
The 75-hole board is shown in Figure~\ref{fig4}.

\noindent
and the d1-complement problem solves by
d3-d1, d5-d3, d7-d5, f7-d7, e5-e7, e3-e5, e1-e3, e8-e6-e4-e2, g6-e6,
c8-e8, e9-e7-e5, i6-g6, i8-i6, j6-h6-f6, g8-g6-e6-e4, k8-i8-g8-e8,
e11-e9-e7, g9-e9, g10-e10-e8-e6, c1-e1-e3-e5-e7, k7-i7-g7, h10-h8,
d10-d8-d6-d4-d2, b6-d6, j9-h9-h7-f7-d7-d5, c4-c6, c7-c5, c2-c4-c6,
z6-b6-d6-d4, x6-z6, y8-y6-a6, z8-z6-b6, a8-a6-c6, b8-b6-d6, d12-d10,
b11-d11-d9, b10-d10-d8, b9-d9-d7-d5-d3-d1.

\noindent
This board has no symmetry whatever, and we have not investigated the
solvability of problems other than the d1-complement.
It appears to us that the square-symmetrical 141-hole board obtained
by attaching 5-arms to the sides of a $9\times 9$ square is not solvable
at the mid-end of the arm, but the 285-hole board obtained by doing
the same to a $15\times 15$ square is solvable everywhere.

\noindent
A 5-arm is the limit.
A board with a 6-arm is unsolvable at the mid-end of the arm whatever
the size and shape of the rest of the board.
The proof is in two stages:
(a) identifying every combination of moves which refills d1 and
clears the rest of the arm, and (b) showing that each leaves a
deficit when measured by the ``golden ratio" resource count developed
by Conway to resolve the problem of the Solitaire Army
(see \cite[chapter 12]{Beasley}).

\section{Five new problems for solution}

Table~\ref{table3} shows the symbols used to describe which holes
are required to be full (a peg is present) at the start and finish
of each problem.
The same symbols are used in \textit{The Ins and Outs} \cite{Beasley}.
A ``marked peg" is one specifically identified,
and generally not allowed to jump until near the end,
when it sweeps all remaining pegs off the board.

\begin{table}[htbp]
\begin{center} 
\begin{tabular}{ | c | l | l |}
\hline
Symbol & Start & Finish \\		
\hline
\hline
(none) & Empty & Empty \\
\includegraphics{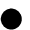} & Full & Empty \\
\includegraphics{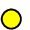} & Marked & Empty \\
\hline
\includegraphics{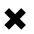} & Empty & Full \\
\includegraphics{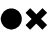} & Full & Full \\
\includegraphics{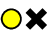} & Marked & Full \\
\hline
\end{tabular}
\caption{Symbols used to describe peg solitaire problems.} 
\label{table3}
\end{center} 
\end{table}

\begin{figure}[htb]
\centering
\epsfig{file=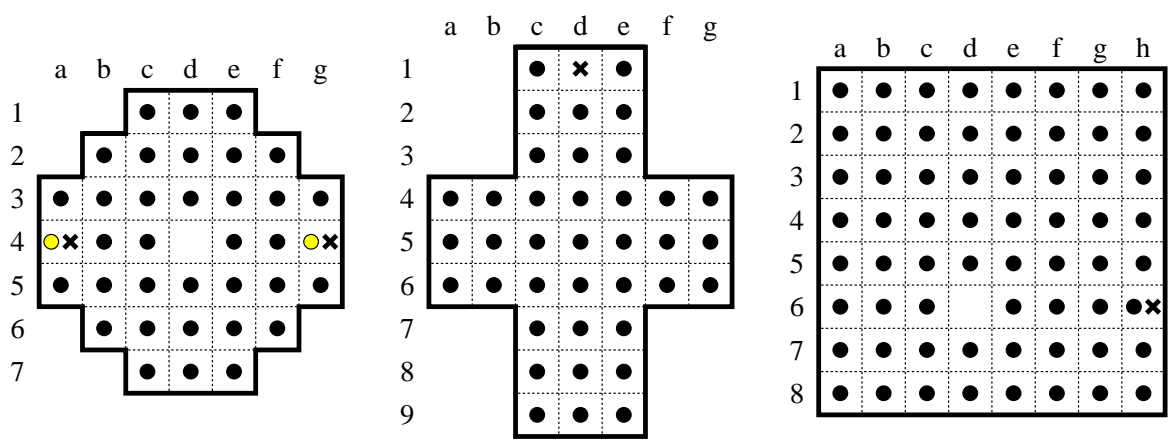}
\caption{(a) Problem 1 on the ``French" board.
(b) Problem 2 on the ``semi-Wiegleb" board.
(c) Problem 3 on the $8\times 8$ square board.}
\label{fig5}
\end{figure}

\noindent
\textbf{Problem 1}\newline
On the 37-hole board, possibly by myself
[John Beasley]\footnote{I am reluctant to make an unqualified claim to this,
because ``vacate d4, finish at a4 and g4" is a natural problem to
try on the 37-hole board and it must have occurred to somebody
to see if it could be done interchanging the pegs originally in these holes,
but I haven't seen it anywhere else.}:
Vacate d4, mark the pegs at a4 and g4, and play to interchange
these pegs and clear the rest of the board (Figure~\ref{fig5}a).

\noindent
\textbf{Problem 2}\newline
On the 39-hole ``semi-Wiegleb" 3-3-3-7-7-7-3-3-3 board,
by George Bell: Vacate d1, and play to finish there
(Figure~\ref{fig5}b).

\noindent
This was discovered in the course of the investigation
described in Section~4.1.
George's computer originally threw out a solution in 24 moves,
my solution by hand took 23;
a subsequent analysis by computer to find the shortest
possible solution got the number down to 21.

\noindent
\textbf{Problem 3}\newline
On the $8\times 8$ board, by John Harris, 1986:
Vacate d6, play to finish at h6 (Figure~\ref{fig5}c).

\noindent
``Here is something I found with 63 poker chips and a chessboard."
John does it in 25 moves, only one more than the number immediately
established as necessary by the $8\times 8$
version of Merson's ``region" analysis.

\begin{figure}[htb]
\centering
\epsfig{file=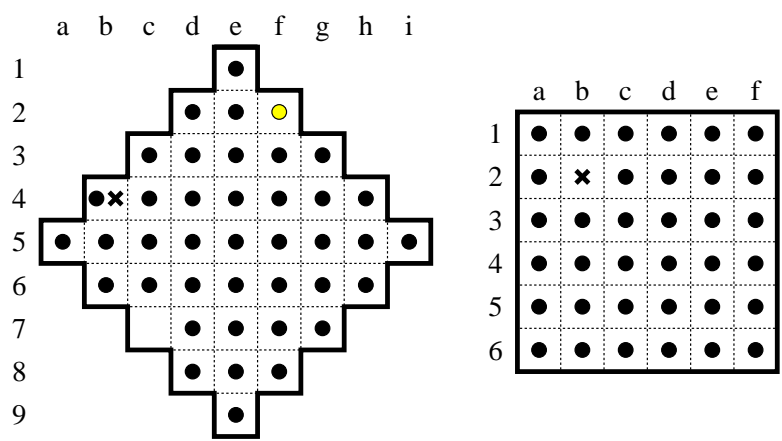}
\caption{(a) Problem 4 on the 41-hole diamond board.
(b) Problem 5 on the $6\times 6$ square board.}
\label{fig6}
\end{figure}

\noindent
\textbf{Problem 4}\newline
On the 41-hole diamond board, by John Harris, 1985:
Allowing diagonal jumps, vacate c7, mark f2,
and play to finish at b4 with a 23-sweep (Figure~\ref{fig6}a).

\noindent
``Can the 41 cell board be cleared in less than 12 moves? Probably.
Is a longer sweep possible on this board?
Don't know, it is possible to set up a 26 peg sweep,
but not if you start with a single vacancy."

\noindent
\textbf{Problem 5}\newline
On the $6\times 6$ board, by John Harris, 1985:
Allowing diagonal jumps, start and finish at b2,
solving the problem in 13 moves and ending with a symmetrical 16-sweep
(Figure~\ref{fig6}b).

\noindent
John's proof that 13 moves are required:
each of the 12 Merson regions around the edge requires a first escape,
and the first jump has to be by a centre peg.
``It is so simple, maybe even a computer could do it!
There could be a 16 peg sweep, 12 move game by starting with the vacancy
somewhere else, but it is unlikely to be symmetrical."

\noindent
Readers are requested to try to solve the problems
for themselves.
This is the best way to gain a full understanding of
any problem.
Problems 1, 4 and 5 are best solved indirectly---first
try to determine the board position before the final sweep(s).
Then, start from the complement of this
board position and  attempt to reduce the board to one peg
at the location of the stating vacancy
(see the ``time-reversal trick" \cite[p. 817--8]{WinningWays}).

\noindent
George Bell has created an
\href{http://www.geocities.com/gibell.geo/pegsolitaire/Tools/GPJ28/GPGPegSolitaireTool.htm}
{interactive JavaScript puzzle} \cite{Bellweb}
where you can try all five problems.

\vskip 30pt 
\begin{footnotesize}

\end{footnotesize}

\newpage
\section*{Solutions}

\begin{figure}[htb]
\centering
\epsfig{file=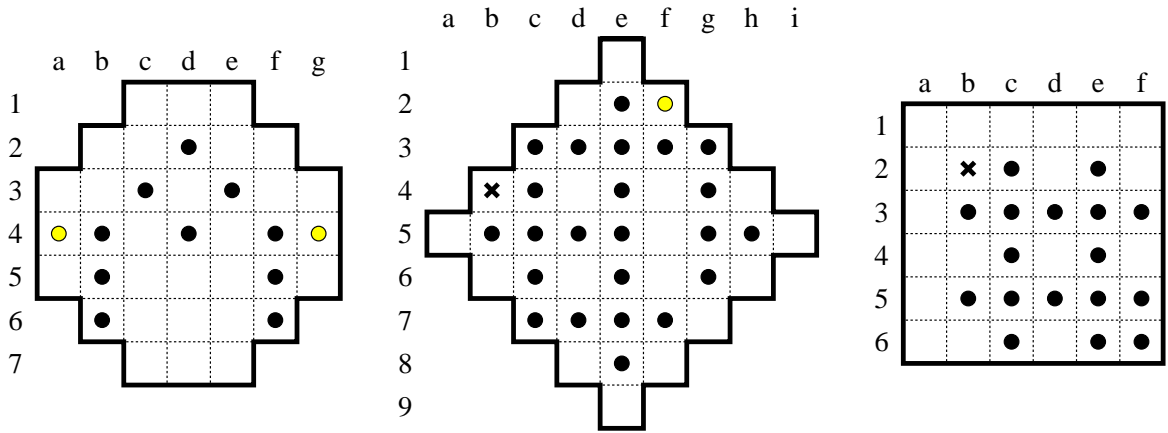}
\caption{(a,b,c) The final sweep positions for Problems 1, 4 and 5.}
\label{fig7}
\end{figure}

\begin{small}
\textbf{Solution to Problem 1} (Figure~\ref{fig5}a)\newline
d2-d4, b2-d2, d1-d3, f2-d2, c4-c2, e3-c3, c2-c4, a3-c3, c4-c2, c1-c3, g3-e3, 
e4-e2, e1-e3, c6-c4, a5-c5, c4-c6, c7-c5, d5-b5, e6-e4, g5-e5, e4-e6, e7-e5, 
d7-d5-f5, after which the board position of Figure~\ref{fig7}a is reached,
and the rest is easy.
\end{small}

\begin{small}
\noindent
\textbf{Solution to Problem 2} (Figure~\ref{fig5}b)\newline
I originally played d3-d1, d5-d3, b5-d5, c3-c5, c1-c3, c6-c4-c2, e5-c5, 
a6-c6, d6-b6, a4-a6-c6-c4, e1-c1-c3-c5, c8-c6-c4, b4-d4-d2, e7-e5, e9-e7, 
e4-e6-e8, g6-e6, d8-d6-f6, g4-g6-e6, c9-e9-e7-e5, e2-e4, f4-d4, f5-d5-d3-d1. 
This was the result of a detailed analysis of debts and surpluses using pencil 
and paper, and had George not told me that the problem was solvable I would have 
assumed it wasn't; indeed, at one point I was sure I had proved it. George's 
computer subsequently reduced the number of moves to 21 by playing d3-d1, d5-d3, 
f4-d4-d2, b5-d5, e6-e4, e3-e5 (6), c7-c5, c9-c7, b4-d4, e1-e3, c2-c4-c6-c8 (11), 
a6-c6, g6-e6, d6-f6, d8-d6-b6 (15), a4-a6-c6, e8-e6-e4-e2, e9-c9-c7-c5-e5 (18), 
g4-g6-e6-e4, c1-e1-e3-e5, f5-d5-d3-d1. 

\noindent
\textbf{Solution to Problem 3} (Figure~\ref{fig5}c)\newline
f6-d6, c6-e6, f8-f6-d6, c8-c6-e6, a8-c8 (5), d8-b8, h8-f8-d8-d6-f6, g6-e6-e8, 
g4-g6-g8, a6-a8-c8 (10), e4-g4, h4-f4, c4-e4-g4, d2-d4, a4-a6-c6-c4-e4-e6 (15), 
b3-d3, c1-c3, a2-a4-c4-c2, a1-c1, d1-b1-b3 (20), f2-f4-f6-d6-d4-d2, f1-d1-d3-f3, 
h2-h4-f4-f2, h1-f1-f3-h3, h6-h8-f8-d8-b8-b6-b4-b2-d2-f2-h2-h4-h6.
Move 8  (g6-e6-e8) is the one that is not an initial exit from one
of the Merson regions. 

\noindent
\textbf{Solution to Problem 4} (Figure~\ref{fig6}a)\newline
e5-c7, c3-e5, f4-d6, g7-e5, h4-f6-d4 (5), c7-e5-g5, e9-c7, b6-d8, 
i5-g7-e9-c7, a5-c3-e5 (10), e1-c3, and we are set up for the sweep
in Figure~\ref{fig7}b,
f2-h4-f4-f2-d2-f4-d4-d2-b4-b6-d4-d6-b6-d8-d6-f6-d8-f8-f6-h6-h4-f6-d4-b4.

\noindent
\textbf{Solution to Problem 5} (Figure~\ref{fig6}b)\newline
d4-b2, a1-c3 (2), b6-d4-b2, a3-a1-c3, a5-a3-c5, d6-b6, a6-c6 (7), f2-d4-b2, 
c1-a1-c3, e1-c1-e3, f4-f2, f1-f3 (12) and we are set up for the sweep
in Figure~\ref{fig7}c,
f6-d6-b6-b4-d6-d4-b4-b2-d4-f6-f4-f2-d2-f4-d4-d2-b2. John uses a binumeric 
notation in order to bring out the symmetry. 
\end{small}

\end{document}